\documentclass{amsart}
\usepackage{amsmath,amssymb}
\usepackage{pstricks}

\newcommand{\todo}[1]{\vspace{5 mm}\par \noindent
\marginpar{\textsc{ToDo}}
\framebox{\begin{minipage}[c]{0.95 \textwidth}
#1 \end{minipage}}\vspace{5 mm}\par}
\renewcommand{\todo}[1]{}
\newcommand{\idiot}[1]{\vspace{5 mm}\par \noindent
\marginpar{\textsc{Note}}
\framebox{\begin{minipage}[c]{0.95 \textwidth}
\tt #1 \end{minipage}}\vspace{5 mm}\par}
\renewcommand{\idiot}[1]{}

\theoremstyle{definition}
\newtheorem{thm}{Theorem}[section]
\newtheorem{thm*}[thm]{Theorem}

\newtheorem{cor}[thm]{Corollary}
\newtheorem{example}[thm]{Example}
\newtheorem{defn}[thm]{Definition}
\newtheorem{lemma}[thm]{Lemma}
\newcommand{\lex}{{\text{Lex}}}

\newcommand{\field}{\Bbbk}

\title{Gotzmann Edge Ideals}
\author{Andrew H. Hoefel}
\begin{document}
\begin{abstract}
Let $P = \field[x_1, \ldots, x_n]$ be the polynomial ring in $n$ variables. 
A homogeneous ideal $I \subseteq P$ generated in degree $d$ is called Gotzmann if 
it has the smallest possible Hilbert function out of all homogeneous ideals
with the same dimension in degree $d$. The edge ideal of a simple graph $G$ on vertices $x_1, \ldots, x_n$ is
the quadratic square-free monomial ideal generated by all $x_ix_j$ where $\{x_i,x_j\}$ is an edge of $G$.
The only edge ideals that are Gotzmann are those edge ideals corresponding to star graphs.
\end{abstract}
\let\thefootnote\relax\footnotetext{The author's research is supported by Killam and NSERC postgraduate scholarships.}
\let\thefootnote\relax
\let\thefootnote\relax\footnotetext{Compiled on \today.}

\maketitle

\section{Introduction}

Throughout this paper, $P = \field[x_1,\ldots, x_n]$ will denote the polynomial 
ring in $n$ variables over the field $\field$ and equipped with the standard grading. 
The \emph{Hilbert function} of a graded $P$-module $M$ is defined to 
be $H(M,d) = \dim_\field M_d$ where $M_d$ is the $d$-th homogeneous component of $M$.

\begin{defn}[Gotzmann Ideal] \label{gotzideal}
A homogeneous ideal $I \subseteq P$ generated in degree $d$ is called \emph{Gotzmann} if 
for all homogeneous ideals $J \subseteq P$ with $H(I, d) = H(J,d)$, we have
\[ H(I,d+1) \leq H(J, d+1).
\]
For brevity we refer to Gotzmann ideals generated in degree $d$ as \emph{$d$-Gotzmann} ideals.
\end{defn}

\idiot{
It is a consequence of the Gotzmann persistence theorem \cite{gotzmann} that out of all homogeneous ideals
with a fixed dimension in degree $d$, Gotzmann ideals have the smallest Hilbert function.
}

Gotzmann ideals are interesting because their Hilbert functions are as small as possible in the following sense.
\begin{thm}
Let $I$ be a $d$-Gotzmann ideal in $P$. For all homogeneous ideals $J \subseteq P$ with $H(I,d) = H(J,d)$, we have
\[	H(I,k) \leq H(J,k) \qquad \text{for all $k$.}
\]
In other words, a $d$-Gotzmann ideal has the smallest Hilbert function out of all homogeneous ideals with the same dimension in degree $d$.
\begin{proof}
This is a well known consequence of the Gotzmann persistence theorem \cite{gotzmann} combined with Macaulay's theorem on Hilbert functions (Theorem \ref{macaulay}).
\end{proof}
\end{thm}

\begin{defn}[Edge Ideal]
Let $G = (V,E)$ be a simple graph on vertices $V\nolinebreak=\nolinebreak\{x_1, \ldots, x_n\}$ and edges $E$. 
The \emph{edge ideal} of $G$ is defined to be 
\[ I(G) = ( x_i x_j \mid \{x_i, x_j\} \in E) \subset P. 
\]
\end{defn}

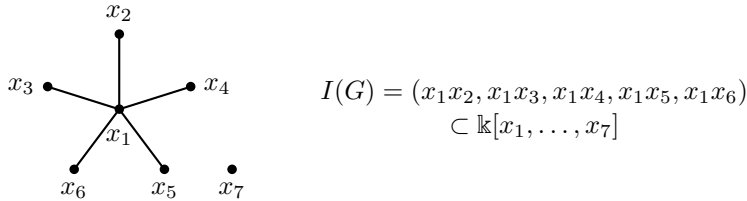
\begin{figure*}[h]
\label{starfigure}
\begin{center}
\begin{pspicture}(-1,-1.25)(7,1.25)
\psdots(0,0)(-.95,0.3)(.95,0.3)(0.6,-.8)(-0.6,-.8)(0,1)(1.5,-.8)
\psline(0,0)(0,1)
\psline(0,0)(-.95,0.3)
\psline(0,0)(.95,0.3)
\psline(0,0)(.6,-0.8)
\psline(0,0)(-.6,-0.8)
\uput[d](0,-0.1){$x_1$}
\uput[u](0,1){$x_2$}
\uput[l](-.95,0.3){$x_3$}
\uput[r](.95,0.3){$x_4$}
\uput[d](0.6,-.8){$x_5$}
\uput[d](-0.6,-.8){$x_6$}
\uput[d](1.5,-.8){$x_7$}
\uput[r](2.5,0){\shortstack{$I(G) = (x_1x_2,x_1x_3, x_1x_4,x_1x_5, x_1x_6)$\\ $ \subset \field[x_1, \ldots, x_7]$}}
\end{pspicture} \\
\caption{A star graph and its edge ideal.}
\end{center}
\end{figure*}

A graph $G$ is called a \emph{star} if there exists a vertex $x_{i_0} \in G$ for which the degree of $x_{i_0}$ is equal to the number of edges in $G$.
See Figure \ref{starfigure} for an example.

The main result of this paper, Theorem \ref{mainResult}, is that the only Gotzmann edge ideals are given by star graphs.

A \emph{lexicographic segment} is a set of monomials of a fixed degree that appear before all 
other monomials of the same degree using the lexicographic order.
Ideals generated by lexicographic segments are fundamental examples of Gotzmann ideals.

Murai and Hibi \cite{gotzmannIdeals} have proven an interesting characterization of all Gotzmann ideals generated by at most $n$ 
homogeneous polynomials. It is easy to see from their result that edge ideals of star graphs are in fact Gotzmann,
but it does not immediately follow from their result that these are the only Gotzmann edge ideals. It is interesting that
to prove that Gotzmann edge ideals only come from star graphs, we will need to consider graphs with more than $n$ edges 
as a separate case (cf. Theorem \ref{boundOnEdgesThm}).

For more information on Gotzmann persistence and the role of lexicographic ideals, we refer the reader to Bruns and Herzog \cite{brunsHerzog} or
Green's lecture on generic initial ideals \cite{genericInitialIdeals}.
Also, Villarreal's book \cite{villarreal} provides a good introduction to edge ideals.

\section{Bounds on Hilbert Functions}

Macaulay's theorem, given below as Theorem \ref{macaulay}, characterizes all possible Hilbert functions of homogeneous ideals in $P$ in terms of
Macaulay representations and pseudo-powers. As a corollary, one can tell if an ideal is $d$-Gotzmann by computing its 
Hilbert function in degrees $d$ and $d+1$.

\begin{defn}[Macaulay Representations and Pseudo-Powers] Let $a\geq 0$ and $d>0$ be integers. Then $a$ can be expressed uniquely as
\[ a = \binom{b_d}{d} + \binom{b_{d-1}}{d-1} + \cdots + \binom{b_1}{1}
\]
where the $b_i$ are integers satisfying $b_d > b_{d-1} > \cdots > b_1 \geq 0$. Note that each $b_i$ depends on both $a$ and $d$. 
This expression for $a$ is called the \emph{$d$-th Macaulay representation} of $a$ and the $b_i$ are called the \emph{$d$-th Macaulay coefficients}.

The \emph{$d$-th Macaulay pseudo-power} and the \emph{$d$-th Kruskal-Katona pseudo-power} of $a$ are defined to be
\[ a^{\langle d \rangle} = \binom{b_d+1}{d+1} + \binom{b_{d-1}+1}{d} + \cdots + \binom{b_1+1}{2}
\]
and
\[ a^{(d)} = \binom{b_d}{d+1} + \binom{b_{d-1}}{d} + \cdots + \binom{b_{1}}{2}
\]
respectively. Here and throughout this paper we assume $\binom{p}{q} = 0$ for $p < q$.
\end{defn}

The Hilbert function of a homogeneous ideal $I \subseteq P$ determines the Hilbert function of its quotient ring $P/I$
in the following way:
\begin{equation} \label{quotienteq}
 H(P/I, d) = H(P,d) - H(I,d) = \binom{n+d-1}{n-1} - H(I,d).
\end{equation}
Consequently, $I$ is a $d$-Gotzmann ideal if and only if $H(P/I,d+1) \geq H(P/J, d+1)$ for all homogeneous ideals $J$ with
$H(P/I, d) = H(P/J,d)$.

\begin{thm*}[Macaulay \cite{macaulay}]
\label{macaulay}
Let $h : \mathbb N \to \mathbb N$ be a function. There is a homogeneous ideal $I \subset P$
with $H(P/I, d) = h(d)$ for all $d \geq 0$ if and only if $h(0) = 1$, $h(1) \leq n$ and $h(d+1) \leq h(d)^{\langle d \rangle}$ for $d \geq 1$.
\end{thm*}

An elegant proof of Macaulay's theorem is given by Green in \cite{restrictionToHyperplanes} and is reproduced in greater context in \cite{brunsHerzog}, \cite{kreuzerRobbiano} and \cite{genericInitialIdeals}.

The following characterization of Gotzmann ideals is a well known corollary of Macaulay's theorem. 
It is often used as the definition of a Gotzmann ideal.

\begin{cor}\label{macaulaycor}
A homogeneous ideal $I\subset P$ generated in degree $d>0$ is Gotzmann if and only if $H(P/I, d+1) = H(P/I, d)^{\langle d  \rangle}$.
\begin{proof}
Let $I$ be $d$-Gotzmann and let $h: \mathbb N \to \mathbb N$ be given by
\[
h(i) = \begin{cases}
	\binom{n+i-1}{n-1} & 0 \leq i < d, \\
	H(P/I,d) & i = d, \\
	h(i-1)^{\langle i-1 \rangle} & i > d.
	\end{cases}
\]
Note that $h(i) = H(P,i)$ for $i < d$ and 
\[ h(d) = H(P/I,d) \leq H(P,d) =  h(d-1)^{\langle d-1 \rangle}.
\]
Since $h$ satisfies the conditions of Macaulay's theorem, there exists an ideal $J$ with $H(P/J, i) = h(i)$ for all $i \geq 0$.
Thus 
\[ H(P/I, d)^{\langle d \rangle } = H(P/J, d+1) \leq H(P/I, d+1)
\] 
as $I$ is Gotzmann. The opposite direction of Macaulay's theorem gives 
\[ H(P/I,d+1) \leq H(P/I,d)^{\langle d \rangle}
\]
and hence we have equality.

The other direction follows immediately from Macaulay's theorem.
\end{proof}
\end{cor}

\idiot{ shortened proof:
\begin{proof}
Let $I$ be $d$-Gotzmann. By Macaulay's theorem there is an ideal $J$ with $H(P/J,d) = H(P/I, d)$ and $H(P/J,d+1) = H(P/J,d)^{\langle d \rangle}$.
Thus $H(P/I, d)^{\langle d \rangle } = H(P/J, d+1) \leq H(P/I, d+1)$. Using the opposite direction of Macaulay's theorem gives
\[ H(P/I,d+1) \leq H(P/I,d)^{\langle d \rangle}
\]
and hence we have equality.

The other direction follows immediately from Macaulay's theorem.
\end{proof}
}

A \emph{simplicial complex} $\Delta$ on ground set $X = \{x_1, \ldots, x_n\}$ 
is any set of subsets of $X$ with the property that $F \in \Delta$ and $E \subset F$ imply $E \in \Delta$. For the purposes of this paper we also assume that $\Delta \neq \emptyset$. 

Sets $F \in \Delta$ are called \emph{faces} and their \emph{dimensions} are given by
$\dim F = |F| - 1$. The \emph{dimension of a complex} is the largest dimension of its faces. 
The \emph{$f$-vector} of a $d$-dimensional simplicial complex is the vector $(f_0, f_1, \ldots, f_d) \in \mathbb N^{d+1}$ where $f_i$
is the number of faces in $\Delta$ with dimension $i$.

Every square-free monomial ideal $I$ corresponds to a simplicial complex $\Delta(I)$ called the 
\emph{Stanley-Reisner complex} of $I$ which is defined as
\[ \Delta(I) = \left\{ F \subseteq \{x_1,\ldots, x_n\} \;\middle\vert\; \prod_{x_i \in F} x_i \notin I \right\}.
\]
The ideal $I$ can be recovered from $\Delta(I)$ by taking the ideal generated by square-free monomials $x_{i_1} x_{i_2} \cdots x_{i_k}$ 
for which $\{x_{i_1}, \ldots, x_{i_k}\}$ is not a face $\Delta(I)$. It is easy to check that this gives a bijective correspondence between 
square-free monomial ideals and simplicial complexes.

The \emph{Stanley-Reisner ring} of a simplicial complex $\Delta$, denoted $\field[\Delta]$, is the quotient ring $P/I$ where 
$\Delta = \Delta(I)$. Stanley-Reisner rings and their Hilbert functions are discussed at length in \cite{brunsHerzog}.

\begin{example}
The Stanley-Reisner complex of the ideal $I = (x_1x_2x_3, x_1x_4)$ contained in $\field[x_1, \ldots, x_4]$ 
has faces $\{x_2,x_3,x_4\}$, $\{x_1,x_2\}$, $\{x_1,x_3\}$ and all subsets thereof. It is two dimensional and its $f$-vector is $(4,5,1)$.
\end{example}

One advantage of Stanley-Reisner complexes is that the Hilbert function of $\field[\Delta]$ can easily be described in terms of
the $f$-vector of $\Delta$. On the other hand, it is often difficult to explicitly describe $\Delta(I)$ when $I$ is a complicated ideal.

The following is Theorem 5.1.7 of \cite{brunsHerzog}.

\begin{thm} \label{bhthm}
Let $\field[\Delta]$ be a Stanley-Reisner ring and let $(f_0, \ldots, f_{\dim \Delta})$ be the $f$-vector of $\Delta$.
Then the Hilbert function of $\field[\Delta]$ is
\[ H(\field[\Delta], d) = \begin{cases}
				1 & d = 0\\
				\sum_{i=0}^{\dim \Delta}f_i\binom{d-1}{i} & d > 0.
			  \end{cases}
\]
\end{thm}

The following theorem places bounds on the growth of the $f$-vector of a simplicial complex. In light of the previous
theorem, bounds on the $f$-vector provide bounds on the Hilbert function of the Stanley-Reisner ring.

\begin{thm*}[Kruskal-Katona \cite{kruskal,katona}]
The vector $(f_0, \ldots, f_d) \in \mathbb N^{d+1}$ is the $f$-vector of a simplicial complex
if and only if 
\[ 0 < f_{k+1} \leq f_k^{(k+1)} \qquad \text{for } 0 \leq k < d.
\]
\end{thm*}

We now show that any Gotzmann square-free monomial ideal achieves the Kruskal-Katona 
bound on its $f$-vector in addition to the Macaulay bound on its Hilbert function.

\begin{thm} \label{sqfreegotz}
Let $I$ be a square-free monomial ideal generated in degree $d$ and let $\Delta(I)$ be its Stanley-Reisner complex. 
Let $(f_0, \ldots, f_{\dim \Delta})$ be the $f$-vector of $\Delta(I)$.
If $I$ is Gotzmann then $f_d = f_{d-1}^{(d)}$.
\begin{proof}
Consider the vector $(f_0, \ldots, f_{d-1}, f_{d-1}^{(d)}) \in \mathbb Z^{d+1}$. As this satisfies the inequalities of the Kruskal-Katona theorem, there
must exist a simplicial complex $\Delta'$ with this $f$-vector. In particular, the Stanley-Reisner ring of $\Delta'$ satisfies
\[	H(\field[\Delta'],d) = \sum_{i=0}^{d-1} f_i\binom{d-1}{i}  = H(P/I, d)
\]
using Theorem \ref{bhthm}.
As $I$ is Gotzmann, $H(P/I, d+1) \geq H(\field[\Delta'],d+1)$. However, applying Theorem \ref{bhthm} to this inequality gives
\[
	\sum_{i=0}^{d} f_i\binom{d}{i} 
	\geq \sum_{i=0}^{d-1} f_i\binom{d}{i}  + f_{d-1}^{(d)} 
\]
and hence $f_d \geq f_{d-1}^{(d)}$. The opposite inequality is given by the Kruskal-Katona theorem
and therefore we have the equality $f_d = f_{d-1}^{(d)}$.
\end{proof}
\end{thm}

\todo{Pick one or none of the following:}

\todo{
In \cite{merminCompressed}, Mermin rephrases the Kruskal-Katona theorem as an analogue of Macaulay's theorem for the ring $Q = P/(x_1^2, \ldots, x_n^2)$.
Similarly, Gotzmann ideals of $Q$ can be defined by replacing $P$ with $Q$ in Definition \ref{gotzideal}. It follows easily from the Kruskal-Katona
theorem that $I \subseteq Q$ is Gotzmann if and only if the Hilbert function of $I$ achieves the bound in the Kruskal-Katona theorem. A consequence
of these definitions and Theorem \ref{sqfreegotz} is that if $I \subset P$ is a $d$-Gotzmann square-free monomial ideal then the image of $I$ in $Q$ is also $d$-Gotzmann.
}

Homogeneous ideals in the quotient ring $Q=P/(x_1^2, \ldots, x_n^2)$, called the \emph{Kruskal-Katona ring}, behave much like their counterparts in $P$.
The Hilbert function of a monomial ideal $I \subseteq Q$ counts the number of square-free monomials contained in each degree of $I^c$, the contraction of $I$ in $P$. If $I$ is a square-free monomial ideal in $P$ then the Hilbert function $h: \mathbb N \to \mathbb N$ of $I +(x_1^2, \ldots, x_n^2) \subseteq Q$ 
is equal to the $f$-vector of $\Delta(I)$ shifted by one index. More precisely, $f_i = h(i+1)$ for $i \geq 0$ where we take $f_i = 0$ for $i > \dim \Delta(I)$.

In \cite{merminCompressed}, Mermin rephrases the Kruskal-Katona theorem as an analogue of Macaulay's theorem for the Kruskal-Katona ring. 
That is, a function $h: \mathbb N \to \mathbb N$ is the Hilbert function of a homogeneous quotient ring of $Q$ if and only if $h(k+1) \leq h(k)^{(k)}$ for all $k$. Gotzmann ideals of $Q$ are defined by replacing $P$ with $Q$ in Definition \ref{gotzideal} and, as in Corollary \ref{macaulaycor}, they are simply those ideals which meet the Kruskal-Katona bound.

After all of these definitions, Theorem \ref{sqfreegotz} can be stated more easily:
\begin{cor}
If $I \subset P$ is a $d$-Gotzmann square-free monomial ideal then the image of $I$ in $Q$ is also $d$-Gotzmann.
\end{cor}

\section{Gotzmann Edge Ideals}

Before proving that only star graphs produce Gotzmann edge ideals, it is first shown in Theorem \ref{boundOnEdgesThm} that Gotzmann edge ideals must have
fewer than $n$ edges where $n$ is the number of vertices in the graph. 
The next two lemmas compute $H(I, 3)$ while assuming that $H(I,2) < n$ 
and that $I$ is Gotzmann or an edge ideal. Finally, in Theorem \ref{mainResult}, these lemmas are used along with the Gotzmann criterion in Corollary \ref{macaulaycor} to prove the main result.

\begin{thm} \label{boundOnEdgesThm}
Let $I = I(G)$ be the edge ideal of a graph $G$ on $n$ vertices with $e$ edges. 
If $I$ is Gotzmann then $e < n$.
\begin{proof}
Let $\Delta  = \Delta(I)$ be the Stanley-Reisner complex of $I$ and let $(f_0,\ldots, f_{\dim \Delta})$
be its $f$-vector. From Theorem \ref{sqfreegotz} we know that $f_2 = f_1^{(2)}$ and so, 
\begin{align*} 
H(P/I,2) &= f_0 + f_1 \qquad\qquad \text{and} \\
H(P/I,3) &= f_0 + 2f_1 + f_1^{(2)}
\end{align*}
from Theorem \ref{bhthm}.

\idiot{
From the Kruskal-Katona theorem we know $f_2 \leq f_1^{(2)}$. 
We now show that if $I$ is Gotzmann then $f_2 = f_1^{(2)}$.

From Theorem \ref{bhthm},
\begin{align*} 
H(P/I,2) &= f_0 + f_1 \qquad\qquad \text{and} \\
H(P/I,3) &= f_0 + 2f_1 + f_2.
\end{align*}

Recursively define a sequence of integers by $s_1 = f_1$ and $s_i = (s_{i-1})^{(i)}$ for $i > 1$. From the definition of the Kruskal-Katona pseudo-power,
one can see that there is an integer $k$ with $s_i = 0$ for $i > k$. 
Thus, the Kruskal-Katona theorem guarantees the existence of a simplicial complex 
$\Delta'$ with $f$-vector $(f_0,f_1, s_2, s_3, \ldots,s_k)$.  The Hilbert function of the Stanley-Reisner ring of $\Delta'$ satisfies
\begin{align*}  
 H(\field[\Delta'],2) &= f_0 + f_1 \quad\qquad\qquad \text{and} \\
 H(\field[\Delta'],3) &= f_0 + 2f_1 + f_1^{(2)}.
\end{align*}

As $I$ is Gotzmann we have
\[ f_0 + 2f_1 + f_1^{(2)}  = H(\field[\Delta'],3) \leq H(P/I, 3) = f_0 + 2f_1 + f_2
\]
and hence $f_1^{(2)} = f_2$ as desired. 
}

By Corollary \ref{macaulaycor}, 
\[ H(P/I, 3) = H(P/I, 2)^{\langle 2 \rangle}
\]
and so,
\begin{equation}\label{eq1}
  f_0 + 2f_1 + f_1^{(2)} = (f_0 +f_1)^{\langle 2 \rangle}. 
\end{equation}

Decompose $f_1$ and $f_0 + f_1$ into their second Macaulay representations as 
\begin{align*}
f_1 = \binom{a_2}{2} + \binom{a_1}{1} \qquad
f_0 + f_1= \binom{b_2}{2} + \binom{b_1}{1}
\end{align*}
where $a_2 > a_1 \geq 0$ and $b_2 > b_1 \geq 0$. Substituting these Macaulay representations
into equation \eqref{eq1} gives
\[
\binom{b_2}{2} + \binom{b_1}{1} + \binom{a_2}{2} + \binom{a_1}{1} + \binom{a_2}{3} + \binom{a_1}{2}
= \binom{b_2+1}{3} + \binom{b_1+1}{2} 
\]
which rearranges and simplifies to
\[ \binom{b_2}{3} + \binom{b_1}{2} 
 = \binom{a_2+1}{3} + \binom{a_1+1}{2}
\]
using the binomial identity $\binom{i}{j} + \binom{i}{j+1} = \binom{i+1}{j+1}$.

These are third Macaulay representations and by the uniqueness of Macaulay representations we have
\[ b_2 = a_2+1 \qquad \text{and} \qquad  b_1 = a_1 +1.
\]

As $I$ is generated in degree two, $f_0 = n$ and so
\begin{align*}
   n &= (f_0+f_1) - f_1 \\
     &= \binom{a_2+1}{2} + \binom{a_1+1}{1}  -  \binom{a_2}{2} - \binom{a_1}{1}\\
     &= a_2 + 1
\end{align*}
again using the binomial identity mentioned earlier. Rearranging gives $a_2 = n-1$ and so $f_1 = \binom{n-1}{2} + \binom{a_1}{1}$. 

The relationship between $e$ and $f_1$ gives
\[ e = \binom{n}{2} - f_1 = \binom{n}{2} -  \binom{n-1}{2} - \binom{a_1}{1} = n-1 - a_1
\]
and since $a_1 \geq 0$ we have $e < n$.
\end{proof}
\end{thm}

\begin{lemma}\label{lexHilbertFunctionLemma}
Let $I \subset P = \field[x_1, \ldots, x_n]$ be a homogeneous ideal generated in degree two and
let $m~=~H(I,2)$. If $m \leq n$ then $I$ is Gotzmann if and only if 
\[ H(I,3) = mn + \frac{1}{2}m - \frac{1}{2} m^2.
\]
\begin{proof}
In the case where $H(I,2)=0$, we have $I=(0)$ and the result clearly holds.

If $H(I,2) \neq 0$ then $H(P/I,2) < H(P,2) = \binom{n+1}{2}$. 
On the other hand, applying equation \eqref{quotienteq} to $P/I$ in degree two gives 
a lower bound as follows:
\[ H(P/I, 2) = H(P,2) - H(I,2) \geq \binom{n+1}{2} - n = \binom{n}{2}.
\]

Thus $H(P/I,2)$ can be written in its second Macaulay representation as 
\begin{equation}\label{mrepeq}
 H(P/I, 2) = \binom{n}{2} + \binom{a}{1} 
\end{equation}  
for some integer $a$ with $n > a \geq 0$.

We compute $a$ by using equation \eqref{quotienteq} once more in degree two which gives
\begin{align*}
 H(I,2) 
   &= H(P,2) - H(P/I,2) \\
   &= \binom{n+1}{2} - \binom{n}{2} - \binom{a}{1} \\
   &= n - a
\end{align*}
and hence $a = n-H(I,2) = n-m$. Replacing $a$ with $n-m$ in equation \eqref{mrepeq} gives
\[ H(P/I,2) = \binom{n}{2} + \binom{n-m}{1}.
\]

By Corollary \ref{macaulaycor}, $I$ is Gotzmann if and only if 
\[ H(P/I,3) = H(P/I,2)^{\langle 2 \rangle} =  \binom{n+1}{3} + \binom{n-m+1}{2}.
\]
Applying \eqref{quotienteq} one last time yields an equivalent condition on $H(I,3)$. Namely, $I$ is Gotzmann if and only if 
\begin{align*} 
H(I,3) 
&= H(P,3) - H(P/I,3)\\
&= \binom{n+2}{3} - \binom{n+1}{3} - \binom{n-m+1}{2}  \\
&= mn + \frac{1}{2}m - \frac{1}{2} m^2.
\end{align*}
\end{proof}
\end{lemma}

Given a graph $G=(V,E)$, a set $S \subseteq V$ is said to be \emph{independent}
if there are no edges $\{u,v\} \in E$ with $\{u,v\} \subseteq S$. Subsets of $V$ which are not independent are called \emph{dependent}.
The faces of the Stanley-Reisner complex $\Delta$ of an edge ideal $I(G)$ are simply the independent sets of $G$. Also note that
the Stanley-Reisner ring $\field[\Delta]$ is equal to the quotient ring $P/I(G)$.

Consequently, the entries of the $f$-vector of $\Delta$ count the number of independent sets of $G$ of a given size.
In particular, $f_1$ is the number of independent sets of size two or, put differently, 
\[ f_1 = \binom{n}{2} - |E(G)|
\]
is the number of non-edges of $G$.

\begin{lemma}\label{edgeIdealHilbertFunctionLemma} Let $G$ be a graph with $e$ edges and $t$ dependent sets of size three. Then 
\[ H(I(G),3) = 2e + t.
\]
\begin{proof}
The monomial basis of $I(G)_3$, the degree three component of $I(G)$, can be partitioned into monomials of the form $x_i^3$, $x_i^2 x_j$ and $x_i x_j x_k$ where $i,j$ and $k$ are distinct.
There are no monomials in $I(G)_3$ of the first type as $I(G)$ is generated by square-free monomials. There are two monomials of type $x_i^2 x_j$ in $I(G)$ for
each edge of $G$ and there is one monomial of type $x_i x_j x_k$ in $I(G)$ for each dependent set of size three.
\end{proof}
\end{lemma}

\begin{thm} \label{mainResult} Let $G$ be a graph. The edge ideal $I(G)$ is Gotzmann if and only if $G$
is a star.
\begin{proof}
Let $G$ be a graph on $n$ vertices and $e$ edges.

We begin by assuming that $I(G)$ is Gotzmann. We know from Theorem \ref{boundOnEdgesThm} that $e < n$.
By Lemma \ref{edgeIdealHilbertFunctionLemma}, $H(I(G),3) = 2e + t$ where $t$ is the number
of dependent sets in $G$ of size three. 

If $G$ contains no edges then $G$ is a star and we are done. So we may assume $e > 0$.

Pick an arbitrary vertex $v$ of $G$ with at least one neighbour. Let $d=\deg v$ be the degree of $v$ and let $H = G \setminus v$ be the graph
obtained from $G$ by deleting $v$. 
Let $e_H = e - d$ denote the number of edges in $H$ and define $t_H$ to be the number
of dependent sets in $H$ of size three.

Let $t_v$ be the number of dependent sets in $G$ of size three which contain $v$.  
We can partition the dependent sets of $G$ into those that contain $v$ and those that are dependent
sets of $H$. Thus, $t = t_v + t_H$. 
\begin{align} 
H(I(G),3) 
&= 2e + t \label{idealhf3} \\
&= 2d + 2e_H + t_H + t_v. \notag
\end{align}

Using the construction in the proof of Corollary \ref{macaulaycor}, Macaulay's theorem guarantees the existence of a 2-Gotzmann ideal $J$ in a polynomial ring in $n-1$ variables with $H(J, 2) = e_H = H(I(H),2)$. 
Thus, 
\begin{align}
2e_H + t_H 
&= H(I(H), 3) \label{graphheq} \\
&\geq H(J, 3) \notag \\
&= (e-d)(n-1) + \frac{1}{2}(e-d)(1-e+d) \notag
\end{align}
using Lemma \ref{edgeIdealHilbertFunctionLemma}, that $J$ is Gotzmann and Lemma \ref{lexHilbertFunctionLemma}.

We now compute $t_v$ -- the number of dependent sets in $G$ of size three containing $v$.
Partition these dependent sets into those that contain two neighbours, one neighbour and no neighbours of $v$.
Every choice of two neighbours of $v$, along with $v$ itself, is dependent. Every choice of a neighbour and a non-neighbour
of $v$, along with $v$, is also dependent. Finally, every choice of two non-neighbours of $v$ which have an edge between them
gives a dependent set of size three when $v$ is included. Thus,
\begin{equation}
 t_v = \binom{d}{2} + d(n-d-1) + |E(G \setminus N(v))| \label{tveq}
\end{equation}
where $N(v)$ is the set of neighbours of $v$.

\idiot{
Taking equation \eqref{idealhf3} and substituting in equations \eqref{graphheq} and \eqref{tveq} for $2e_H + t_H$ and $t_v$ respectively
gives the following after much simplification:
\[ H(I(G),3) = \left(ne + \frac{1}{2}e - \frac{1}{2}e^2\right) +  (d-1)(e-d)  + |E(G \setminus N(v))|.
\]
}

Taking equation \eqref{idealhf3} and substituting in inequality \eqref{graphheq} and equation \eqref{tveq} for $2e_H + t_H$ and $t_v$ respectively
gives the following:
\begin{align*}
H(I(G),3) &= 2d + (2e_H + t_H) + t_v \\
          &\geq 2d + (e-d)(n-1) + \frac{1}{2}(e-d)(1-e+d) \\ 
	  &\phantom{{}\geq 2d} {}+ \binom{d}{2} + d(n-d-1) + |E(G \setminus N(v))| \\
	  &= ne - \frac{1}{2}e^2 + de -d^2 -\frac{1}{2}e +d + |E(G \setminus N(v))| \\
	  &= \left(ne + \frac{1}{2}e - \frac{1}{2}e^2\right) +  (d-1)(e-d)  + |E(G \setminus N(v))|.
\end{align*}

\idiot{
\begin{align*}
&= 2d + (e-d)(n-1) + \frac{1}{2}(e-d)(1-e+d) + |E(G \setminus N(v))| + d(n-d-1) +\frac{1}{2}d(d-1) \\
&= -(e-d) + 2d + n(e-d) + \frac{1}{2}(e-d)(1-e+d) + |E(G \setminus N(v))| + d(n-d-1) +\frac{1}{2}d(d-1) \\
&= -(e-d) +de -d^2 + 2d + n(e-d) + \frac{1}{2}(e -e^2 -d -d^2) + |E(G \setminus N(v))| + d(n-1) +\frac{1}{2}d(d-1) \\
&= (d-1)(e-d) + 2d + n(e-d) + \frac{1}{2}(e -e^2 -d -d^2) + |E(G \setminus N(v))| + d(n-1) +\frac{1}{2}d(d-1) \\
&= (d-1)(e-d) + 2d + ne + \frac{1}{2}e(1-e) - nd + \frac{1}{2}(-d -d^2) + |E(G \setminus N(v))| + d(n-1) +\frac{1}{2}d(d-1) \\
&= (d-1)(e-d) + ne + \frac{1}{2}e(1-e)  + |E(G \setminus N(v))| \\
&= (d-1)(e-d) + H(\lex_e^2,3)  + |E(G \setminus N(v))| 
\end{align*}
}

As $I(G)$ is 2-Gotzmann, $H(I(G),3) = ne + \frac{1}{2}e - \frac{1}{2}e^2$ from Lemma \ref{lexHilbertFunctionLemma} and 
so we have 
\[  (d-1)(e-d)  + |E(G \setminus N(v))| \leq 0
\]
However, both terms above are non-negative and hence we have 
\begin{align*}
   (d-1)(e-d) &= 0 \qquad\text{and}  \\
   |E(G \setminus N(v))| &= 0.
\end{align*}
These equations hold for every choice of $v \in V(G)$ with degree $d \geq 1$ and so, 
$d$ is either $0$, $1$ or $e$ for every vertex of $G$.
If every vertex $v$ in $G$ has degree $0$ or $1$ and $e \neq 1$ then $G$ has more than one connected component and so $|E(G \setminus N(v))|$ cannot be zero. Thus, $G$ must have some vertex $v$ with degree $d = e$ and hence $G$ is a star.

Conversely, if $G$ is a star then every dependent set of $G$ must contain the common vertex $v$ of all edges.
That is $t = t_v = \binom{e}{2} + e(n-e-1) + |E(G \setminus N(v))|$ from equation $\eqref{tveq}$. However, $G\setminus N(v)$ contains no edges.
Therefore, 
\begin{align*}
H(I(G),3) 
&= 2e + t \\
&= 2e + \binom{e}{2} + e(n-e-1) \\
&= ne + \frac{1}{2}e - \frac{1}{2}e^2
\end{align*}
and hence $I(G)$ is Gotzmann by Lemma \ref{lexHilbertFunctionLemma}.  
\end{proof}
\end{thm}

The author would like to thank Sara Faridi and Emma Connon for their careful reading of this manuscript.

\bibliographystyle{plain}
\bibliography{gotzmann_edge_ideals}

\end{document}